\documentclass[12pt,letterpaper]{article}
\usepackage{amsmath}
\usepackage{amsfonts}
\usepackage{amssymb}
\usepackage{latexsym}

\newenvironment{itemizz}{\begin{itemize}\setlength{\itemsep}{-1mm}} %
{\end{itemize}}                              

\newcommand\eop{$\ \ {\vcenter
   {\hrule
   \hbox{\vrule height 9pt \kern 9pt \vrule height 9pt}
   \hrule}}$\vskip 1.0 pt}
\newenvironment{proof}{{\bf Proof.}}{\eop\medskip}
\newenvironment{proofof}[1]{\medskip \textbf{Proof of #1.}}{\eop\medskip}

\newtheorem{theorem}{Theorem}[section]
\newtheorem{definition}[theorem]{Definition}
\newtheorem{lemma}[theorem]{Lemma}
\newtheorem{corollary}[theorem]{Corollary}
\newtheorem{proposition}[theorem]{Proposition}

\newcommand\CC{{\mathcal C}}

\newcommand\duX{{\widehat{X}}}
\newcommand\duN{{\widehat{N}}}
\newcommand\duG{{\widehat{G}}}

\newcommand\TTT{{\mathbb T}}
\newcommand\T{{\mathbb T}}
\newcommand\ZZZ{{\mathbb Z}}
\newcommand\Z{{\mathbb Z}}
\newcommand\RRR{{\mathbb R}}

\newcommand\Q{{\mathbb Q}}
\newcommand{\g}{\mathfrak{g}}

\newcommand\Hom{{\mathrm{Hom}}}
\newcommand\q{\mbox{\rm qc}}
\newcommand{\gcl}{\mathfrak{g}\mathrm{cl}}

\newcommand{\res}{\mathord{\restriction}}  % less space around it.
\newcommand\bohr{\mathsf{b}}  % Bohr compactification
\newcommand\sbl[1]{\langle#1\rangle}   % subloop
\newcommand\hull[1]{\langle#1\rangle}   % subloop
\newcommand{\un}{\underline}

\begin{document}

\title{Characterizing Subgroups of Compact Abelian Groups}
\author{Dikran Dikranjan\thanks
{Partially supported by INDAM (GNSAGA)}
\ and Kenneth Kunen\thanks
{Partially supported by NSF Grant DMS-0097881.}}

\maketitle

\vskip -21pt
\vbox{}

\begin{abstract}
We prove that every countable subgroup of a compact metrizable 
abel\-ian group has a characterizing set. As an application, we answer
several questions on maximally almost periodic (MAP) groups
 and give a characterization of the class of (necessarily MAP)
abelian topological groups whose Bohr topology has countable pseudocharacter.   
\end{abstract}

\section{Introduction}
 
We shall write our abelian groups additively, and 
we view the circle $\TTT$ as $\RRR/\ZZZ$;
then, for $z \in \TTT$, $\|z\|$ is the distance to the nearest integer,
so $0 \le \|z\| \le 1/2$.  All topological groups are assumed to 
be Hausdorff unless otherwise noted.

For a topological abelian group $X$,
$\duX$ denotes the group of all continuous characters on $X$
(i.e., the continuous homomorphisms from $X$ to $\T$),
and $\duX$ is given the compact open topology. 
When $X$ is locally compact, 
the \emph{Pontryagin Duality Theorem} lets us identify
$X$ with $\widehat{\duX\,}$; see \cite{FOL,HR, RUD}.
For most of the results in this paper, $X$ will be compact,
so that $\duX$ will be discrete.

One can use countable sequences or countable sets of characters
to define subgroups of $X$ as follows:

\begin{definition} Let $X$ be a topological abelian group:
\begin{itemizz}
  \item[a.]    \textup{\cite{DMT}}  For a sequence
$\un{u}=\langle u_n : n \in \omega \rangle$ of 
elements of $\duX$, let $s_{\un{u}}(X)$
be the set of all $x \in X$ such that $u_n(x)\to 0$.
If $H = s_{\un{u}}(X)$, we say that $\un{u}$ \emph{characterizes} $H$.
  \item[b.] \textup{\cite{HK1,HK2}} For a countably 
infinite subset $B$ of $\duX$, let
 $\CC_B(X)$ be the set of all $x \in X$ such that
$\langle \varphi(x) : \varphi \in B\rangle$ converges to $0$ in $\T$.
If $H = \CC_B(X)$, we say that $B$ \emph{characterizes} $H$.
\end{itemizz}
\end{definition}

Note that $s_{\un{u}}(X)$ and $\CC_B(X)$ are both subgroups of $X$.
We study here which subgroups of $X$ can be characterized as
an $s_{\un{u}}(X)$ or a $\CC_B(X)$.
Observe that
$s_{\un{u}}(X)$ and $\CC_B(X)$ are minor variants of the same notion:

\begin{lemma}
\label{lemma-sB}
Let $X$ be any compact topological abelian group.
\begin{itemizz}
\item[1.] $C_B(X) = s_{\un{u}}(X)$ for any $\un u$ which is
a 1-1 enumeration of $B$.
\item[2.] If $H$ is closed
and $|X : H|$ is finite, then $H = s_{\un{u}}(X)$ for some $u$.
\item[3.] Every $\CC_B(X)$ is a Haar null set.
\item[4.] If $|X : H|$ is infinite and
$H = s_{\un{u}}(X)$, then $H = \CC_B(X)$ for some $B$.
\end{itemizz}
\end{lemma}
\begin{proof}
(1) is clear from the definitions.
For (2), $H^\perp = \{\varphi \in \duX : \varphi(H) = \{0\}\}$
is finite (of size $|X : H|$),
and $H = s_{\un{u}}(X)$ for any $\un u : \omega \to H^\perp $
which lists each character
in $H^\perp $ infinitely often.
(3) follows easily from the orthogonality of the characters;
see, e.g., \cite{BDMW,CTW,HK1,HK2}.  (4) will be proved 
in Section \ref{sec-proofs}.
\end{proof}

Note also that every $C_B(X)$ and $s_{\un{u}}(X)$ is a Borel
set (in fact, an $F_{\sigma\delta}$), and every Borel subgroup
of finite index is closed and open and of positive Haar measure.
By Lemma \ref{lemma-sB}, characterization as an $s_{\un{u}}(X)$
is equivalent to characterization as a $\CC_B(X)$ except in
the trivial case of subgroups of finite index.

The subgroups $s_{\un{u}}(X)$ have been studied by many authors,
especially in the case $X=\T$, where
$\duX$ is identified with $\ZZZ$
(see \cite{A,BDMW2,BDMW3,BDS,BS,DdS,E1,KL,Larcher}).
When $x \in \TTT$ is a non-torsion element, 
the behavior of the sequences of the form 
$\langle u_n(x) : n \in \omega \rangle$ is related to
Diophantine approximation and dynamical systems
(more specifically, the Sturmian flow \cite{PS}), as well as to the 
study of precompact group topologies with converging sequences
\cite{BDMW,BDMW4,DMT}. 
 
Some instances of characterization in $\TTT$
were described  in Armacost's book \cite{A}, although this book
predates the ``characterization'' terminology.  Specifically,
for each prime $p$, the
Pr\"ufer subgroup $\Z_{p^\infty}$ of $\T$ is characterized
by $\{p^n : n \in \omega\}$.
Also, let $D = \{n! : n \in \omega\}$.
Then the set of rational points of $\T$ (i.e., $\Q/\Z$) is properly
contained in $\CC_D$ (in fact, this $\CC_D$ is uncountable).
The elements of $\CC_D$ are called the
\emph{topologically torsion elements} of $\T$ in \cite{A};
they are described further in \cite{D, DdS,DPS}. 

It is easy to see that $\Q/\Z$ can be characterized by some set of characters:

\begin{proposition}
\label{prop-facts}
If $B = \{k\cdot n! : 0 <  k \le n \in \omega\}$, then
$\Q/\Z = \CC_B(\T)$.
\end{proposition}

It is not always so easy to write down a characterizing set
explicitly, but we shall show:

\begin{theorem}
\label{aaa} 
Every countably infinite subgroup of a compact 
metrizable abelian group $X$ has a characterizing set.
\end{theorem}

This theorem resolves Problem 5.3 from \cite{DMT}.
The proposition and theorem will be proved in Section \ref{sec-proofs}.
A number of important special cases of the theorem are
already in the literature.  In particular, when  $X = \TTT$, 
it was proved by
B\' \i r\' o,   Deshouillers, and S\' os \cite{BDS}.
Earlier results by Larcher \cite{Larcher}
and Kraaikamp and Liardet \cite{KL} relate the characterizability
of infinite cyclic subgroups
$\hull{\alpha}$ of $\T$ to the continued fraction
representation of $\alpha$; see also \cite{BDMW2}.
It was proved in \cite{DMT} that every cyclic subgroup of the
group $\widehat{\Z_{p^\infty}}$ 
of $p$-adic integers has a characterizing set.
Furthermore, B\' \i r\' o \cite{B2} has proved
Theorem \ref{aaa} in the case that the subgroup is dense in $X$,
finitely generated, and torsion-free.
Also, Beiglb\"ock, B\'\i r\'o, S\'os, and Winkler (unpublished) have
proved Theorem \ref{aaa} by a different method.

Note that in Theorem \ref{aaa}, $X$ must be metrizable 
(equivalently, second countable),
since $\CC_B(X)$ always contains $\bigcap_{\varphi \in B} \ker(\varphi)$,
and this set will be uncountable if $X$ is not metrizable.
We do not know a simple general criterion for deciding whether
$H \le X$ is characterizable.  
For closed $H$, this is easy; $H$ is of the form
$s_{\un{u}}(X)$ iff $H$ is a $G_\delta$
(see Proposition \ref{prop-closed} below),
and then,
applying Lemma \ref{lemma-sB},
$H$ is of the form $\CC_B(X)$ iff $H$ is a $G_\delta$
of infinite index in $X$.
In Section \ref{sec-proofs}, we shall prove the following
theorem, which shows that not all $F_\sigma$
subgroups can be characterized, even in metrizable compact groups:

\begin{theorem}
\label{thm-strict}
Suppose that $X$ is a compact abelian group
and $H = \bigcup_n F_n$, where each $F_n \le X$ is closed
and each $F_n \lneqq F_{n+1}$.
Then the following are equivalent:
\begin{itemizz}
\item[a.] $H = \CC_B(X)$ for some countably infinite $B \subseteq \duX$.
\item[b.] $H = s_{\un{u}}(X)$ for some $\un u : \omega \to \duX$.
\item[c.] For some $m$: $X/F_m$ is metrizable and $| F_{n+1} : F_n |$
is finite for all $n \ge m$.
\end{itemizz}
\end{theorem}

Next, we give an application of
Theorem \ref{aaa} to non-compact groups.
If $(G, \tau)$ is a topological abelian group, we let
$(G, \tau^+)$ be its \emph{Bohr modification}.  
So, $\tau^+$ is the coarsest topology which makes all
(continuous) characters of $(G, \tau)$ continuous.
Clearly, $\tau^+$ is coarser than $\tau$.
When $\tau$ is clear from context, we refer to $(G, \tau^+)$
as $G^+$.  It is easy to find Hausdorff $G$ for which $G^+$ is indiscrete,
but we are primarily interested in the case where $G^+$ is
also Hausdorff; such $G$ are called
\emph{MAP} (\emph{maximally almost periodic}).

A topological space $X$ has \emph{countable pseudocharacter}
iff every singleton $\{x\}$ is a $G_\delta$ set.
This implies that $X$ is a $T_1$ space, which, in the case of
topological groups, is equivalent to being Hausdorff.

\begin{corollary}
\label{b}
For a MAP topological abelian group $G$ the following are equivalent: 
\begin{itemizz}
\item[a.] \emph{Every} countable $H \le G$ is of the form $s_{\un{u}}(G)$
for some $\un u : \omega \to \duG$.
\item[b.] \emph{Some} countable $H \le G$ is of the form $s_{\un{u}}(G)$
for some $\un u : \omega \to \duG$.
\item[c.] $G^+$ has countable pseudocharacter.
\end{itemizz}
\end{corollary}
\begin{proof}
$(a) \to (b)$ is obvious.  For $(b) \to (c)$, note that
$K = \bigcap_n \ker(u_n)$ is a $G_\delta$ set in $G^+$,
and $K \le H$, so $K$ is countable. 
Since $G^+$ is Hausdorff, it follows that $\{0\}$ is
a  $G_\delta$ set in $G^+$.

For $(c) \to (a)$:  Since $\{0\}$ is a  $G_\delta$ set in $G^+$,
there is a countable $\Phi \subseteq \duG$ such that
$\{0\} = \bigcap\{\ker(\varphi) : \varphi \in \Phi\}$.
Let $\Delta : G \to \T^\Phi$ be the evaluation map:
$\Delta(x)(\varphi) = \varphi(x)$.
Then $\Delta$ is 1-1 and continuous.
Applying Theorem \ref{aaa} and Lemma \ref{lemma-sB}.1,
let $\un{v}=\langle v_n : n \in \omega \rangle$ be a
sequence of characters of $\T^\Phi$  such that
$\Delta(H) = s_{\un{v}}(\T^\Phi)$.  If $u_n = v_n \circ \Delta$,
then $H = s_{\un{u}}(G)$.
\end{proof}

Note that if one did not assume that $G$ is MAP in Corollary \ref{b},
either of condition (a) (just using $\{0\}$ = $s_{\un{u}}(G)$)
or condition (c) would imply that $G$ is MAP anyway.

Since having a characterizing sequence (or set)
is a rather restrictive property,
one can relax it in the following way, following \cite{DMT}:

\begin{definition}
For a subgroup $H$ of a topological abelian group $X$,
\[
\gcl(H) = \bigcap\left\{s_{\un{u}}(X) :
\ \ u : \omega \to \duX  \ \&\  H  \subseteq s_{\un{u}}(X)\right\} \ \ .
\]
$H$ is \emph{$\g$-closed} iff $\gcl(H) = H$.
\end{definition}

Of course, $H$ is $\g$-closed whenever $H = s_{\un{u}}(X)$ for some 
$\un u $, but a $\g$-closed subgroup need not be of the
form $s_{\un{u}}(X)$.  For example, 

\begin{proposition}
\label{prop-closed}
Suppose that $H$ is a closed subgroup of the compact $X$.  Then:
\begin{itemizz}
\item[1.] $H$ is $\g$-closed.
\item[2.] $H = s_{\un{u}}(X)$ for some $\un u$ iff 
$H$ is a $G_\delta$ set.
\end{itemizz}
\end{proposition}
\begin{proof}
For (1), fix $x \notin H$, and then fix $\varphi \in \duX$ with
$\varphi(x) \ne 0$ but $\varphi(H) = \{0\}$.  If $u_n = \varphi$ for
all $n$, then $x \notin s_{\un{u}}(X)$.  Thus, $x \notin \gcl(H)$.
For $(2)\leftarrow$, if $H$ is a closed $G_\delta$, then
$H^\perp = \{\varphi \in \duX : \varphi(H) = \{0\}\}$
is countable, and then $H = s_{\un{u}}(X)$ for any $\un u : \omega \to H^\perp $
which lists each character in $H^\perp $ infinitely often.
For $(2)\rightarrow$, let $K = \bigcap_n \ker(u_n)$,
and let $\pi: X \to X/K$ be the natural quotient map.
Then $X/K$ is second countable, and $\pi(H)$ is closed and hence
a $G_\delta$ in $X/K$, so $H = \pi^{-1} \pi(H)$ is 
a $G_\delta$  in $X$.
\end{proof}

Likewise, Theorem \ref{aaa} holds only for compact metrizable $X$,
but metrizability is not needed for:

\begin{corollary}
\label{cor-aaa} 
Every countably infinite subgroup of a compact 
abelian group $X$ is $\g$-closed.
\end{corollary}
\begin{proof}
Fix a countable $H$ and $x \notin H$.
Then, let $\Gamma$ be a countable subgroup of $\duX$ which separates
the elements of $H \cup \{x\}$, let
$K = \bigcap \{\ker(\psi) : \psi \in \Gamma\}$,
and let $\pi: X \to X/K$ be the natural quotient map.
Then $\pi$ is 1-1 on $H \cup \{x\}$, and one can 
apply Theorem \ref{aaa}  in $X/K$ to fix $B \subseteq \widehat{X/K}$
with $\pi(H) = \CC_B(X/K)$.  If
$\{v_n : n \in \omega\}$ is a 1-1 listing of $B$
and $u_n = v_n \circ \pi$, then 
$H \le s_{\un{u}}(X)$ and $x \notin s_{\un{u}}(X)$.
\end{proof}

This generalizes Theorem 1.2 of \cite{DMT}, which proved the result for
cyclic subgroups, and it answers Problem 5.1 and Question 5.2 of \cite{DMT}.
Note that Corollary \ref{cor-aaa} actually holds for all MAP $X$,
since one can apply the corollary in the Bohr compactification, $\bohr(X)$.
We do not know if every $F_\sigma$ subgroup of a compact $X$
is $\g$-closed, although, in contrast to Theorem \ref{thm-strict},

\begin{proposition}
Suppose that $X$ is a compact abelian group
and $H = \bigcup_n F_n$, where each $F_n \le X$ is closed
and each $F_n \le F_{n+1}$.
Then $H$ is $\g$-closed.
\end{proposition}
\begin{proof}
Fix $x\in X\backslash H$.
Since $F_n$ is closed, we can choose
$u_n \in \duX$ such that $F_n \le \ker (u_n)$ and $\|u_n(x)\| \ge 1/4$.
Then 
$H \le s_{\un{u}}(X)$ and $x \notin s_{\un{u}}(X)$.
\end{proof}

$F_{\sigma\delta}$ subgroups certainly need not be  $\g$-closed.
In fact, whenever $X$ is infinite, there is an
$F_{\sigma\delta}$ subgroup $H$ such that $H$ is a Haar null set
and $\gcl(H) = X$; see \cite{HK2}, which extends
earlier arguments in \cite{BDMW, HK1}.

The next corollary says that every totally bounded countable metrizable 
group topology can be ``encoded'' by
means of a \emph{single} $\tau$-convergent zero-sequence:

\begin{corollary}\label{repe}
Let $(G,\tau)$ be a countable totally bounded metrizable abelian group.
Then there exists a sequence $\un u = \langle u_n : n \in \omega \rangle$
with $u_n \to 0$ in $(G,\tau)$ such that $\tau$ is the finest 
totally bounded group topology on $G$ in which $u_n \to 0$.
\end{corollary}

\begin{proof}
Let $G_d$ denote the group $G$
with the discrete topology, and let $X = \widehat{(G_d)}$.
Let $H \le X$ be the set of $\varphi : G \to \T$ which are
$\tau$-continuous.  Then $H$ is countable and $X$ is metrizable.
Identifying $G_d$ with $\duX$,
apply Theorem \ref{aaa} to fix $u_n \in G$ with 
$H = s_{\un{u}}(X)$.
Now, suppose that $\tau'$ is a strictly finer
totally bounded group topology on $G$.
Then there is a $\tau'$-continuous character $x$ of $G$
which is not $\tau$-continuous, so that
$x \in X$ but $x \notin H$.  Then
$x(u_n) = u_n(x) \not\to 0$ in $\TTT$, so $u_n \not\to 0$ in $(G, \tau')$.
\end{proof}

It is easy to see that countablity is necessary; that is,
if $(G,\tau)$ is a totally bounded metrizable abelian group such that
$\tau$ is the finest totally bounded group topology on $G$ in which
some sequence $u_n \to 0$, then $G$ is countable. 

Corollary \ref{repe} should be compared to the fact that the finest
\emph{Hausdorff} group topology that makes a
given sequence converge to 0 is never 
Fr\'echet-Urysohn \cite{ZP}.  
In the case $G = \ZZZ$, the encoding sequence $\un u$ can
often be displayed explicitly.
For example, if $\tau$ is the $p$-adic topology,
one can take $u_n = p^n$.
If $\tau$ is the profinite topology, then
one can let $\un u$ enumerate
$\{k\cdot n! : 0 \le k \le n \in \omega\}$ by Proposition \ref{prop-facts}.

\section{The Proofs}
\label{sec-proofs}

We begin by relating characterizable subgroups of $N$
to characterizable subgroups of $X$, where $N$ is a subgroup of $X$:

\begin{lemma}
\label{lemma-subgroup}
Suppose that $H \le N \le X$,  where  $X$ is a
compact abelian group,  $N$ is closed in $X$,
$X/N$ is metrizable,
and $H = \CC_B(N)$ for some countably infinite $B \subseteq \duN$.
Then $H = \CC_D(X)$ for some countably infinite $D \subseteq \duX$.
\end{lemma}
\begin{proof}
Let $N^\perp = \{\varphi \in \duX : \varphi(N) = \{0\}\}$,
which is countable.

If $|N^\perp| = \aleph_0 $:
For $\varphi \in B$, choose $\widetilde \varphi \in \duX$
such that $\widetilde \varphi \res N = \varphi$.
Let $D = \{\widetilde \varphi : \varphi\nobreak \in\nobreak B\}  \cup N^\perp $.
Then $\CC_D^X \cap N = \CC_B^N$, and if $y \in X \backslash N$,
then $y \notin \CC_{N^\perp}^X$, so $y \notin \CC_D^X $.

If $|N^\perp| = k < \aleph_0 $:
Then $|X : N| = k$.
Let $D = \{\psi \in \duX : \psi \res N \in B\}$.
For $\varphi \in B$, $|\{\psi \in D : \psi \res N = \varphi\}| = k$,
so  $\CC_D^X \cap N = \CC_B^N$.
Now, fix $y \in X \backslash N$, and fix any $\alpha \in N^\perp$ such
that $\alpha(y) \ne 0$.  Note that $\psi \in D$ iff $\alpha + \psi \in D$,
so $\langle \psi(y) : \psi \in D\rangle$ 
cannot converge to $0$.
\end{proof}

\begin{proofof}{Lemma \ref{lemma-sB}.4}
We have $X$ compact, $H \le X$, 
$|X : H|$ infinite, and $H = s_{\un{u}}(X)$.
We shall produce a countably infinite $D \subseteq \duX$ such that
$H = \CC_D(X)$.
Let $N = \overline H$ and let
$K = \bigcap\{\ker(u_n) : n \in \omega\}$.
Then $K \le H \le N \le X$, and $X / K$ is metrizable.
If $H$ is closed in $X$, then  $X/H$ is metrizable and infinite,
so  $H = \CC_{H^\perp}(X)$.

Now, assume that $H$ is not closed, so $H \ne N$.
Let $v_n = u_n \res N$, so $H = s_{\un{v}}(N)$.
Let $B = \{v_n : n \in \omega\} \subseteq \duN$. 
Then $B$ is infinite (since otherwise $H$ would be closed),
and $\{n :u_n = \psi\}$ is finite for each non-zero $\psi \in \duN$
(since $N = \overline H$), so $H = \CC_B(N)$.
We are now done by Lemma \ref{lemma-subgroup}.
\end{proofof}

\begin{proofof}{Proposition \ref{prop-facts}}
It is clear that $\Q/\Z \subseteq \CC_B(\T)$.
Now, any $x \in [0,1]$ can be written in a factorial expansion
as $x=\sum_{n=1}^\infty c_{n}/(n+1)!$, where the integers $c_n$ 
satisfy $0\leq c_{n}\leq n$.
Note that $\sum_{n=k}^\infty n/(n+1)! = 1/k!$.
Thus, if $x \notin \Q$, then  $c_n\notin \{ 0, n\}$ for infinitely many $n$.
For any $n$, we have $n!\, x \equiv y_n \pmod 1$, where
$y_n = n! \sum_{m=n}^\infty c_{m}/(m+1)!$.  When $c_n\notin \{ 0, n\}$,
\[
\frac{1}{n+1} =
n! \left[ \frac{1}{(n+1)!}\right]
\le y_n \le
n! \left[ \frac{n-1}{(n+1)!} + \frac{1}{(n+1)!}\right]
= \frac{n}{n+1} \ \ ,
\]
so we can find a $k$ with $0 < k \le n$ such that 
$\|k n! x\| = \|k n! y_n\| \ge 1/4$. 
When $x \notin \Q$, there are infinitely many
such $k,n$, so $x \notin \CC_B(\T)$.
\end{proofof}

We now prove three lemmas for Theorem \ref{aaa}.
As is common in metric spaces, $B(x;\varepsilon)$ denotes an open
ball, and $N(F;\varepsilon)$ denotes the open set
$\bigcup\{B(x;\varepsilon) : x \in F \}$.

\begin{lemma}
\label{lemma-metric}
Let $(X; d)$ be a metric space,  $F$ any countably infinite subset of $X$,
and $F_n$ finite subsets such that $F_n \nearrow F$.  
Then there are positive $\varepsilon_n \searrow 0$ such that
for all $x \notin F$, there are infinitely many $n$
with $x \notin N(F_n; \varepsilon_n)$.
\end{lemma}
\begin{proof}
Choose $\varepsilon_n \searrow 0$ so that
$2 \varepsilon_n < d(u,v)$ whenever $u,v$ are distinct elements of 
$F_{n+1}$.
Now, fix $m$ and assume that $x \in N(F_n; \varepsilon_n)$ for all $n \ge m$;
we show that $x \in F$.
For each $n \ge m$, choose $u_n \in F_n$ such that $d(x,u_n) < \varepsilon_n$.
Then $d(u_n, u_{n+1}) \le d(x,u_n) + d(x,u_{n+1}) <
\varepsilon_n + \varepsilon_{n+1} \le 2 \varepsilon_n$,
so $u_n = u_{n+1}$.
Thus, the $u_n$ are all equal to some $u \in F$, and $x = u$.
\end{proof}

If $G$ is a discrete group, then $G$ is dense in its Bohr
compactification; if this fact is stated in terms of the
compact group $X = \duG$, we get the following lemma,
which generalizes an old result of Kronecker;
see  \cite[ p.~1188]{C},
\cite[Cor.~26.16]{HR}, or \cite{HZ}.

\begin{lemma}
\label{lemma-Kron}
Let $X$ be a compact abelian group,
$x_1,\ldots,x_t$ a finite list of elements of $X$,
$\psi:X\to \T$ a possibly discontinuous homomorphism,
and $\delta>0$.   Then there is a $\varphi\in \widehat{X}$ such
that $\|\psi(x_i)-\varphi(x_i)\|<\delta$ for $i=1,\ldots, t$.
\end{lemma}

Using this, one can give a direct proof of Corollary \ref{cor-aaa}, 
without using Theorem \ref{aaa}. 
We have $H < X$, where $H$ is countable and $X$ is compact.
Fix $x \notin H$ and list $H$ as $\{e_n : n \in \omega\}$.
For each $n$, $\;x \notin \sbl{ e_0, e_1, \ldots, e_n }$ implies
that there is a homomorphism $\psi_n : X\to \T$ such that
$\|\psi_n(x)\| > 1/4$ and $\ker(\psi_n)$ contains $e_0, e_1, \ldots, e_n$.
By Lemma \ref{lemma-Kron}, there is a $u_n \in \duX$ such that
$\|u_n(x)\| > 1/4$ and $\|u_n(e_j)\| < 1/n$ for $j = 0,1,\ldots, n$.
Then $H \le s_{\un{u}}(X)$ and $x \notin  s_{\un{u}}(X)$.

\begin{definition}
Let $X$ be any topological abelian group and
$E \subseteq X$.
The \emph{$m$-quasi-convex hull} $\q_m(E)$ of $E$ is the set 
\[
\left\{ x \in X:\ \  \forall \varphi \in \duX \, \left[
\forall e \in E  \, [ \|\varphi(e)\| \leq 2^{-m-2} ]
\ \ \Longrightarrow\ \  \|\varphi(x)\| \leq 1/4 \right ]  \right\}
\]
$\q(E) = \q_0(E)$.
\end{definition}
The \emph{quasi-convex hull} $\q(E)$ is discussed in
\cite{Au, Ban}.
Note that $\q_m(E)$ gets bigger as $m$ gets bigger, and
$E \subseteq \q_0 (E)$.

\begin{lemma}\label{lemma-F} Let $m\in \omega$ and let $E$ be a finite
subset of a compact abelian group $X$. Then: 
\begin{itemizz}
  \item[a.] $\q_m(E) \subseteq \sbl{E}$.
  \item[b.] $\q_m(E)$ is finite.
\end{itemizz}
\end{lemma}

\begin{proof} For (a), fix $x \not \in \sbl{E}$. There exists a possibly 
discontinuous homomorphism $\psi$ such that
$\psi(E)=\{0\}$ and $\|\psi(x)\| > 1/4$.
Apply Lemma \ref{lemma-Kron} to find a $\varphi\in\duX$ with 
$|\varphi(e)| <2^{-m-2}$ for each $e \in E$ and 
$\|\varphi(x)\| > 1/4$, so that $x\notin \q_m(E)$.

For (b), let $H = \sbl{E}$.
Then $H$ is an internal direct product of cyclic groups; 
say $H \cong \sbl{a_1} \times \cdots \times \sbl{a_n}$,
where $a_1, \ldots, a_n \in H$.  For each $j$, the order $o(a_j)$
is either a positive integer or $\infty$.
Fix a positive integer $M \ge \max\{o(a_j) : o(a_j) \ne \infty\}$
such that also each $e \in E$ is of the form $\sum_{j = 1}^n \mu_j a_j$,
where each $|\mu_j| \le M$.
Let $F_0$ be the set of all elements $\sum_{j = 1}^n \nu_j a_j$
such that each $|\nu_j| \le 2^{m+1}M$.  Then $F_0$ is
finite, and we show $\q_m(E) \subseteq F_0$.

Fix $x = \sum_{j = 1}^n \nu_j a_j \notin F_0$, with each 
$|\nu_j| < o(a_j)$.  Then, fix $\ell$ such that $|\nu_\ell| > 2^{m+1}M$.  Then 
$o(a_\ell) > |\nu_\ell| > M$, so $o(a_\ell) = \infty$.
We can then define $\psi \in  \Hom(H, \TTT)$ so that
$\psi(\sum_{j = 1}^n \mu_j a_j) = \mu_\ell / (2 \nu_\ell)\pmod 1\;$
for all $\sum_{j = 1}^n \mu_j a_j \in H$.
Then $\psi(x) = 1/2$, but for each $e  = \sum_{j = 1}^n \mu_j a_j \in E$,
we have $\|\psi(e)\| < M \cdot 1 / (2^{m+2}M) = 2^{-m-2}$. By
Lemma \ref{lemma-Kron},
there is a   $\varphi\in \duX$ such that $\|\varphi(x)\|> 1/4$, but 
$\|\varphi(e)\|< 2^{-m-2}$ for each  $e  \in E$.
Thus, $x \notin \q_m(E)$.
\end{proof}

This lemma actually holds for all MAP $X$,
since one can apply the lemma in  $\bohr(X)$.
For $m = 0$, this lemma was proved in \cite{Au}
(see Lemma 7.10 and Theorem 7.11).

The inclusion $\q_m(E) \subseteq \sbl{E}$ can fail
when $E$ is infinite.  For example, 
$\q_m(E)=\overline{E}$ for every 
$m\in \omega$ whenever $E=\sbl{E} \le X$.

\begin{proofof}{Theorem \ref{aaa}}
Let $X$ be a compact metric abelian group with $G = \duX$.
Let $E$ a countably infinite subgroup
of $X$. In view of Lemma \ref{lemma-subgroup}, we may assume that
$E$ is dense in $X$. Let $E_n \nearrow E$, with each $E_n$ finite.
Let $F_n = \q_n(E_n)$.  Then $E_n \subseteq F_n \subseteq E$,
$F_n \nearrow E$, and $F_n$ is finite by Lemma \ref{lemma-F}.
Let $d$ be a metric for $X$. Applying Lemma \ref{lemma-metric},
choose positive $\varepsilon_n \searrow 0$ such that
for all $x \notin E$, there are infinitely many $n$
with $x \notin N(F_n; \varepsilon_n)$.

Let $A_n=\{\varphi \in \duX:  \forall e \in E_n \,[ \|\varphi(e)\| 
\le 2^{-n-2} ]\}$.
Applying the definition of  $\q_n(E_n)$,
\[
\forall x \notin F_n \; \exists \varphi \in A_n \;
\left[ \|\varphi(x)\| > 1/4\right] \ \ .
\]
Hence,
$\bigcup_{\varphi\in A_n}\{x\in X:\|\varphi(x)\| > 1/4\} \supseteq
X \setminus  N(F_n; \varepsilon_n)$.
Since $X \setminus  N(F_n; \varepsilon_n)$
is compact, we can choose a finite
$B_n \subseteq A_n$ such that 
\[
\forall x \notin N(F_n, \varepsilon_n)
\; \exists \varphi \in B_n \; [ \|\varphi(x)\| > 1/4] \ \ .
\]
Set $B = \bigcup_{n \in \omega} B_n$ and note that 
$B$ must be infinite, since (using $\overline E = X$)
$$
\forall \varphi \ne 0 \; [ | \{n : \varphi \in A_n\} | < \aleph_0 ] \ . \eqno(*)
$$
Clearly, $E \subseteq \CC_B$. To prove that 
$E = \CC_B$ fix $x \notin E$.  Then $x \notin N(F_n; \varepsilon_n)$ for
infinitely many $n$.  For these $n$, choose $\varphi_n \in B_n$
so that $ \|\varphi_n(x)\| > 1/4$.  By $(*)$,
$\{ \varphi_n : x \notin N(F_n; \varepsilon_n)\}$ is infinite,
so $x \notin \CC_B$.
\end{proofof}

Next, we prove three lemmas for Theorem \ref{thm-strict}.

\begin{lemma}
\label{lemma-closed}
Suppose that $X$ is a compact abelian group,
$F$ a closed subgroup, and $B$ a countably infinite subset
of $\duX$.  Then $F \le \CC_B$ iff $F \le \ker(\varphi)$
for all but finitely many $\varphi \in B$.
\end{lemma}
\begin{proof}
For the non-trivial direction, assume that $F \le \CC_B$
but $F \not\le \ker(\varphi)$ for infinitely many $\varphi \in B$.
Since $\CC_B$ gets bigger as $B$ gets smaller, we may assume 
that $\varphi \res F \ne 0$ for all $\varphi \in B$, and
that one of the following two cases holds:

Case I:  The $\varphi\res F = \psi \in \widehat F$
for all $\varphi \in B$:  But then any $x \in F \backslash \ker(\psi)$
will be in $F \backslash \CC_B$, contradicting $F \le \CC_B$.

Case II:  The $\varphi\res F$ for $\varphi \in B$ are all different.
Let $D = \{\varphi\res F : \varphi \in B\}$.  Then
$\CC_B \cap F = \CC_D$, which is a null set in $F$
by Lemma \ref{lemma-sB},  contradicting $F \le \CC_B$.
\end{proof}

\begin{lemma}
\label{lemma-notin}
Suppose that $X$ is a compact metric abelian group, with
an invariant metric $d$.
Let $H \le X$, where $H = \bigcup_n F_n$, each $F_n$ is closed,
and each $F_n < F_{n+1} < X$.
Assume that $y_n \in F_{n+1} \backslash F_n $ are chosen so that
each $d(y_{n+1}, 0) \le (d(y_n, F_n))/3$.
Define $x = \sum_{k=0}^\infty y_k$.  Then $x \notin H$.
\end{lemma}
\begin{proof}
By induction on $k$, we have $d(y_{n+k}, 0) \le (d(y_n, F_n))\cdot 3^{-k}$
for all $n \ge 0$ and $k > 0$.
This shows that the sum defining $x$ really converges, and also
lets us define $x_n = \sum_{k=n}^\infty y_k$; so $x = x_0$.
Now, fix $n$, and we show that $x \notin F_n$.
Since $(x - x_n) \in F_n$, it is sufficient to show that $x_n \notin F_n$.
Now $y_n \notin F_n$ and $x_n = y_n + x_{n+1}$.
Also, $d(x_{n+1}, 0) \le \sum_{k=1}^\infty d(y_{n+k}, 0) \le
d(y_n, F_n) \cdot \sum_{k=1}^\infty 3^{-k} = d(y_n, F_n)/2$.
Thus, $d(x_n, F_n) \ge d(y_n, F_n) - d(x_{n+1}, 0) \ge d(y_n, F_n)/2 >0$.
\end{proof}

Finally, we need the (trivial) converse to Lemma \ref{lemma-subgroup}:

\begin{lemma}
\label{lemma-supergroup}
Suppose that $H \le N \le X$,  where  $X$ is a
compact abelian group,  $N$ is closed in $X$,
$H$ is not closed in $X$, and
$H = \CC_D(X)$ for some countably infinite $D \subseteq \duX$.
Then $H = \CC_B(N)$ for some countably infinite $B \subseteq \duN$.
\end{lemma}
\begin{proof}
By Lemma \ref{lemma-sB}.1, $H = s_{\un{u}}(X)$, where each $u_n \in \duX$.
But then also $H = s_{\un{v}}(N)$, where each $v_n = u_n \res N \in \duN$.
We now get $B$ by Lemma \ref{lemma-sB}.4.
\end{proof}

\begin{proofof}{Theorem \ref{thm-strict}}
$(a) \leftrightarrow (b)$ follows from Lemma \ref{lemma-sB},
since
$|X : H|$ is infinite.
For $(c) \to (a)$, we may,
replacing $X$ by $X/F_m$ assume that $F_m = \{0\}$.
But then $H$ is countable, so the result follows by Theorem \ref{aaa}.

For $(a) \to (c)$: 
Assume $(a)$, and let $K = \bigcap\{\ker(\varphi) : \varphi \in B\}$.
Then $K \le H$ and $X/K$ is metrizable.
Observe that $K \le F_m$ for some $m$; otherwise, we have
$K \cap F_n \nearrow K$,  where each $K \cap F_n \lneqq  K$.
If some $|K : K \cap F_n |$ is finite, this is clearly a contradiction,
while if all $|K : K \cap F_n |$ are infinite, this is a contradiction
by the Baire Category Theorem.

We are now done if we can show that $| F_{n+1} : F_n |$
is finite for all but finitely many $n$; so, assume that this is false
and we shall derive a contradiction.
Re-indexing, we may assume that $| F_{n+1} : F_n |$ infinite for each $n$,
and we may assume that $F_0 = \{0\}$.
We may also assume that $X = \overline H$; if not, replace 
$X$ by $N = \overline H$ and apply Lemma \ref{lemma-supergroup}.
Then $\widehat X = F_0^\perp > F_1^\perp > F_2^\perp \cdots$,
and $\bigcap_n F_n^\perp = \{0\}$.

We may assume that $X$ is metrizable; otherwise replace $X$
by $X / K$.
We may assume that $0 \notin B$,
so that $B = \bigcup_n B_n$,
where $B_n = B \cap F_n^\perp \backslash F_{n+1}^\perp$.
By Lemma \ref{lemma-closed}, each $B_n$ is finite.

We shall now produce an
$x \in \CC_B(X) \backslash H$, contradicting $H = \CC_B(X)$.

Let $d$ be an invariant metric for $X$.
Inductively choose 
$y_n \in F_{n+1} \backslash F_n $ so that
$\|\varphi(y_n)\| \le 2^{-n}$ for all
$\varphi \in B_0 \cup \cdots \cup B_n$ and
each $d(y_{n+1}, 0) \le (d(y_n, F_n))/3$; this is possible
because each $F_n$ is nowhere dense in $F_{n+1}$.
Let $x = \sum_{k=0}^\infty y_k$.  Then $x \notin H$ by Lemma \ref{lemma-notin}.

If $\varphi \in B_n$, then
$\|\varphi(x)\| \le \sum_{k=0}^\infty \|\varphi(y_k)\| =
\sum_{k=n}^\infty \|\varphi(y_k)\| \le
\sum_{k=n}^\infty 2^{-k}  = 2^{1-n}$.
Thus $\langle \varphi(x) : \varphi \in B\rangle$ converges to $0$ in $\T$.
\end{proofof}

\textsc{
Dipartimento di Matematica e Informatica,
Universit\`{a} di Udine,
Via della Scienze 206,
33100 Udine,
Italy}

\textit{Email address}: \verb+dikranja@dimi.uniud.it+

\bigskip

\textsc{Department of Mathematics,
University of Wisconsin, Madison, WI 53706, USA}

\textit{Email address}: \verb+kunen@math.wisc.edu+

\textit{URL:} \verb+http://www.math.wisc.edu/~kunen+


\begin{thebibliography}{99}
  
\bibitem{A} D. L. Armacost, {\it The structure of locally compact Abelian
groups}, Monographs and Textbooks in Pure and Applied Mathematics,
{\bf 68}, Marcel Dekker, Inc., New York, 1981.

\bibitem{Au} L. Aussenhofer,  
{\em Contributions to the duality theory of abelian topological groups 
and to the theory of nuclear groups}, Dissertationes Math. (Rozprawy Mat.) 384 (1999). 

\bibitem{Ban}
W. Banaszczyk, {\em Additive subgroups of topological vector spaces},
 Lecture Notes in Mathematics, 1466. Springer-Verlag, Berlin, 1991.

\bibitem{BDMW} G. Barbieri, D. Dikranjan, C. Milan, H. Weber,
{\it Answer to Raczkowski's questions on convergent sequences}, Top.
Appl. {\bf 132} (2003), 89--101.

\bibitem{BDMW2} G. Barbieri, D. Dikranjan, C. Milan and H. Weber,{\it 
Topological torsion related to some recursive sequences of integers}, 
submitted.

\bibitem{BDMW3}G. Barbieri, D. Dikranjan, C. Milan and H. Weber, {\it 
$\mathfrak t$-dense subgroups of topological Abelian groups}, 
submitted.

\bibitem{BDMW4} G. Barbieri, D. Dikranjan, C. Milan and H. Weber, 
{\it Convergent sequences in precompact group topologies}, submitted.

\bibitem{B1} A. B\' \i r\' o, {\em Characterizing sets for subgroups 
of compact groups I: a special case}, preprint.

\bibitem{B2} A. B\' \i r\' o, {\em Characterizing sets for subgroups 
of compact groups II: general case}, preprint.

\bibitem{BDS} A. B\' \i r\' o, J.-M.  Deshouillers and V. S\' os,
{\em Good approximation and characterization of subgroups of
$\mathbb R/\mathbb Z$},
Studia Sci. Math. Hungar. {\bf 38} (2001), 97--113.

\bibitem{BS} A. B\' \i r\' o and V. S\' os,
{\em  Strong characterizing sequences in simultaneous Diophantine
approximation}, J. Number Theory {\bf 99} (2003), 405--414.
  
\bibitem{C} W. Comfort {\em Topological groups,} in: Handbook
         of Set-Theoretic Topology, edited by K. Kunen and J. E. Vaughan,
          North Holland, Amsterdam $\cdot$ New York $\cdot$
         Oxford (1984), p.1184. 

\bibitem{CTW} W.W. Comfort, F.J. Trigos-Arrieta and T.S. Wu,
{\it The Bohr compactification, modulo a metrizable subgroup},
Fundamenta Math. {\bf 143} (1993), 119-136, Correction: {\bf 152} 
(1997), 97-98.

\bibitem{D} D. Dikranjan, {\it Topologically torsion elements of
topological groups}, Topology Proc., {\bf 26}, 2001-2002, pp.
505--532.

\bibitem{DdS} D. Dikranjan and R. Di Santo, {\it On Armacost's quest on
topologically torsion elements}, Communications Algebra, {\bf 32} 
(2004), 133--146.

\bibitem{DMT} D. Dikranjan, C. Milan and A. Tonolo {\it 
A characterization of the MAP abelian groups}, J. Pure Appl. Algebra (2004) 
to appear.


\bibitem{DPS} D. Dikranjan, Iv. Prodanov and L. Stoyanov, {\it Topological groups. Characters, dualities and minimal group topologies},
 Monographs and Textbooks in Pure and Applied Mathematics 
{\bf 130}, Marcel Dekker, Inc., New York, 1990. 

 
\bibitem{E1} H.G. Eggleston, {\it Sets of fractional dimensions which
occur in some problems of number theory,}
    Proc. London Math. Soc. (2) {\bf 54}, (1952) 42--93.

\bibitem{FOL} G. B.  Folland,
{\it A Course in Abstract Harmonic Analysis},
CRC Press, 1995.

\bibitem{HK1} J. Hart and K. Kunen,
\textit{Limits in function spaces and compact groups},
Topology Appl., to appear.

\bibitem{HK2} J. Hart and K. Kunen, {\em Limits in compact abelian groups}, to appear.
  
\bibitem{HR} E.~Hewitt and K.~Ross, {\it Abstract Harmonic Analysis},
Volume I, Springer-Verlag, Berlin-G\"ottingen-Heidelberg 1979.

\bibitem{HZ} E. Hewitt and H. Zuckerman, {\em A group-theoretic method in approximation theory}, 
Ann. of Math. (2) {\bf 52} (1950), 557--567.

\bibitem{KL} C. Kraaikamp and P. Liardet, {\it
Good approximations and continued fractions},  Proc. Amer. Math.
Soc. {\bf 112} (1991), 303--309.

\bibitem{Larcher} G. Larcher, {\it A convergence problem connected with the
continued fractions}, Proc. Amer. Math. Soc., {\bf 103} (1988), 718-722.

\bibitem{PS} K. Petersen and L. Shapiro, 
{\em Induced flows}, 
Trans. Amer. Math. Soc. {\bf 177} (1973), 375--390.

\bibitem{RUD}  W. Rudin, \textit{Fourier Analysis on Groups},
Interscience Publishers, 1962.

\bibitem{ZP} E. G. Zelenyuk and I. V. Protasov, {\it Topologies
on abelian groups}, Math. USSR Izvestiya {\bf 37} (1991), 445--460.
Russian original: Izvestia Akad. Nauk SSSR {\bf 54} (1990), 1090--1107.
\end{thebibliography}
\end{document}